\documentclass[oneside, reqno] {amsart}

\usepackage{xypic}
\input xy
\xyoption{all}
\usepackage{epsfig}
\usepackage{amsthm}
\usepackage{amssymb}
\usepackage{amsmath}
\usepackage{amscd}
\usepackage{color}
\usepackage[T1]{fontenc}
\usepackage{upgreek}
\usepackage[font=scriptsize]{caption}
\usepackage{wrapfig}

\usepackage{graphicx}
\usepackage{subfigure}

\newcommand{\bg}{\begin{equation}}
\newcommand{\ed}{\end{equation}}
\newcommand{\bga}{\begin{eqnarray}}
\newcommand{\eda}{\end{eqnarray}}
\newcommand{\pf}{\textbf{Proof:\ }}

\def\cbdu{\par{\raggedleft$\Box$\par}}

\newtheorem {Theorem}  {Theorem}

\numberwithin{Theorem}{section}

\newtheorem {Lemma}[Theorem]  {Lemma}
\newtheorem {Proposition}[Theorem]{Proposition}
\theoremstyle{definition}
\newtheorem{Definition}[Theorem]{Definition}
\theoremstyle{remark}

\expandafter\chardef\csname pre amssym.def
at\endcsname=\the\catcode`\@ \catcode`\@=11
\def\undefine#1{\let#1\undefined}
\def\newsymbol#1#2#3#4#5{\let\next@\relax
 \ifnum#2=\@ne\let\next@\msafam@\else
 \ifnum#2=\tw@\let\next@\msbfam@\fi\fi
 \mathchardef#1="#3\next@#4#5}
\def\mathhexbox@#1#2#3{\relax
 \ifmmode\mathpalette{}{\m@th\mathchar"#1#2#3}%
 \else\leavevmode\hbox{$\m@th\mathchar"#1#2#3$}\fi}
\def\hexnumber@#1{\ifcase#1 0\or 1\or 2\or 3\or 4\or 5\or 6\or 7\or 8\or
 9\or A\or B\or C\or D\or E\or F\fi}

\font\teneufm=eufm10 \font\seveneufm=eufm7 \font\fiveeufm=eufm5
\newfam\eufmfam
\textfont\eufmfam=\teneufm \scriptfont\eufmfam=\seveneufm
\scriptscriptfont\eufmfam=\fiveeufm

\catcode`\@=\csname pre amssym.def at\endcsname

\newcounter{remark}
\renewcommand{\div}{\mbox{div}}

\def  \12  {{\frac{1}{2}}}

\def\build#1_#2^#3{\mathrel{\mathop{\kern 0pt#1}\limits_{#2}^{#3}}}

\numberwithin{equation}{section}

\begin{document}

\title[Stationary even active scalar equations]{Non-unique stationary solutions of even active scalar equations}

\author [Mimi Dai]{Mimi Dai}

\address{Department of Mathematics, Statistics and Computer Science, University of Illinois at Chicago, Chicago, IL 60607, USA}
\email{mdai@uic.edu}

\author [Chao Wu]{Chao Wu}

\address{Department of Mathematics, Statistics and Computer Science, University of Illinois at Chicago, Chicago, IL 60607, USA}
\email{cwu206@uic.edu}

\thanks{The authors are partially supported by the NSF grants DMS--2009422 and DMS--2308208.}

\begin{abstract}

We study a class of active scalar equations with even non-local operator in the drift term. Non-trivial stationary weak solutions in the space $C^{0-}$ are constructed using the iterative convex integration approach.

\bigskip


\hspace{0.02cm}CLASSIFICATION CODE: 35Q35, 35Q86, 76D03.
\end{abstract}

\maketitle

\section{Introduction}
\label{sec-int}
We consider the active scalar equation on $\mathbb T^d$
\begin{equation}\label{ase}
\begin{split}
\partial_t\theta+ u \cdot\nabla \theta+\Lambda^\gamma \theta=&\ 0, \\
u=&\ T[\theta],\\
\nabla\cdot u=&\ 0
\end{split}
\end{equation}
where $T$ is a Calder\'on-Zygmund operator with even Fourier symbol. Since $\nabla\cdot T[\theta]= 0$, without loss of generality, we assume 
\[T[\theta]=\nabla\times \Gamma^{-1}[\theta]\]
where $\Gamma^{-1}$ is the inverse operator of $\Gamma$. The operator $\Gamma$ has Fourier symbol $m(k)$ which is odd and homogeneous of degree 1. Such active scalar equations with even drift operator arise from several physical contexts, such as the incompressible porous media (IPM), magneto-geostrophic (MG) model, etc. In these physical examples, $\theta$ denotes an unknown scalar function and $u$ stands for the velocity field. In (\ref{ase}) a dissipation term with generalized Laplacian is present, with $\gamma>0$. 

Classical well-posedness problems have been studied for the aforementioned active scalar equations in many articles, for instance, see \cite{CGO, FV1, Moff} and references therein. The current note concerns the non-uniqueness of weak solutions for (\ref{ase}). It is known that uniqueness of a solution to a nonlinear PDE is only guaranteed if the solution is regular enough. One usually does not expect uniqueness to hold for weak solutions. Although there are different ways to show non-uniqueness, the recent development of convex integration techniques first innovated by Nash \cite{Nash}  in 1950's has shown it is a systematic and generic approach to construct weak solutions for PDEs which violate uniqueness or physical conservation laws. The successful implementation of convex integration schemes in PDEs started from the pioneering work of De Lellis and Sz\'ekelyhidi \cite{DLS1, DLS2} for Euler equations. Since then we have seen a real blossom of the techniques applied to many fluid equations with breakthrough outcomes, including the resolution of the Onsager conjecture for the Euler equations \cite{BDLSV, Is, GR}. For more work in this regard, we refer the reader to the survey paper \cite{BV19}.

Due to the special nonlinear structure of the active scalar equations, we face more obstacles in exploiting the application of convex integration schemes to these equations. One main difficulty is that there is not a suitable family of stationary solutions as building blocks which possesses nice properties as those for the Euler equations. Using the convex integration scheme within the Tartar framework, the author of \cite{Shv} constructed infinitely many bounded weak solutions for (\ref{ase});  non-unique bounded weak solutions for the 2D IPM equation were also obtained in \cite{CFG} and \cite{Sze}; it was shown the existence of non-unique mixing solutions for the 2D IPM in some particular geometry settings in \cite{CCF} and \cite{CFM}.
 Obtaining non-unique weak solutions for \eqref{ase} in H\"older space was contributed in the work \cite{IV}. In particular, the authors developed new ideas to realize the key cancellations in the iterative convex integration process to achieve $C_x^{\alpha}$ regularity with $\alpha<\frac{1}{4d+1}$. 
 
\subsection{Main result} 
In this paper we will construct non-trivial stationary (time independent) weak solutions for (\ref{ase}) on $\mathbb T^d$ for $d\geq 2$. The analysis is presented in 2D, i.e. $d=2$ for simplicity.  Therefore we consider the stationary active scalar equation on $\mathbb T^2$
\begin{equation}\label{S-ase}
\begin{split}
u \cdot\nabla \theta+\Lambda^\gamma \theta=&\ 0, \\
u=T[\theta]=&\ \nabla^{\perp}\Gamma^{-1}[\theta].
\end{split}
\end{equation}

Denote $\mathbb P_{\leq \lambda}$ by the standard Littlewood-Paley projection operator with frequency support below the frequency $\lambda$.

\begin{Definition}\label{def}
A function $\theta:\mathbb T^2\to \mathbb R$ is said to be a stationary weak solution of (\ref{S-ase}) if 
\begin{equation}\notag
-\int_{\mathbb T^2}\theta u\cdot \nabla\psi \, dx+\int_{\mathbb T^2} \theta \Lambda^{\gamma}\psi\, dx=0
\end{equation} 
is satisfied for any function $\psi\in C^\infty(\mathbb T^2)$ satisfying $\psi=\mathbb P_{\leq \Lambda_0}\psi$ for some frequency number $\Lambda_0$.
\end{Definition}

Our main result is:
\begin{Theorem}\label{thm}
Let $0<\gamma<2-\alpha$ and $\alpha<1$. There exists a non-trivial weak solution $\theta$ of (\ref{S-ase}) with $\Lambda^{-1}\theta\in C^{\alpha}(\mathbb T^2)$.
\end{Theorem}

\subsection{Relevant previous work}
In fact, using the method of convex integration, non-uniqueness has been obtained previously for stationary fluid equations, including active scalar equations and Navier-Stokes equations.  For the stationary Navier-Stokes equations on $\mathbb T^d$ with $d\geq 4$, it was shown in \cite{Luo} that there exists a non-trivial weak solution in $L^2$. Non-unique stationary weak solutions were constructed in \cite{CKL} (without external forcing) and \cite{DP} (with external forcing) for the surface quasi-geostrophic (SQG) equation, which is of active scalar equation with an odd drift operator. Comparing stationary equations to their evolutionary cases, it is in general more challenging to construct weak solutions with high regularity, while implementing the convex integration approach. The reason is that the time dimension can be regarded as additional flexibility, since the principal idea of convex integration is to exploit the flexibility of the underlying equations.  

Although it is known from \cite{IV} that non-unique weak solutions with H\"older regularity exist for the evolutionary even active scalar equations without forcing, it remains open to show non-uniqueness of weak solutions in $L^2$ or $C^\alpha$ for $\alpha\geq 0$ for the unforced stationary even active scalar equations. On the other hand, unlike the SQG for which weak solutions can be defined in $H^{-\frac12}$, it seems not applicable to define weak solutions for even active scalar equations with negative regularity index. This motivates the study of a particular type of weak solutions as in Definition \ref{def}.

In the direction of investigating forced equations within the framework of convex integration, the flexibility of allowing a forcing has been revealed in \cite{BHP, BHP2, DF, DP, DP2}. In \cite{BHP} the authors were able to obtain non-unique weak solutions with regularity higher than Onsager's critical $\frac13$ regularity for forced Euler equation. While non-unique weak solutions with sharp regularity were constructed in \cite{BHP2, DP2} for forced SQG. Recently, the authors of \cite{DF} obtained non-unique weak solutions for time dependent forced even active scalar equations in space $C_t^0C_x^\alpha$ for $\alpha<\frac1{2d+1}$, which has higher spatial regularity compared to the solutions obtained in \cite{IV}. However, surprisingly, allowing an external forcing does not help to improve the regularity of convex integration weak solutions for the stationary even active scalar equations. This indicates that the nonlinear structures of the even and odd active scalar equations (and Euler equation) are intrinsically different.

\subsection{Organization of the paper}
We provide a proof of Theorem \ref{thm} in Section \ref{sec-no-force}. 
We begin with a simple convex integration construction for the unforced active scalar equation (\ref{S-ase}), which yields a sequence of approximating solutions $\{\theta_q\}$ with $\Lambda^{-1}\theta_q\in C^{1-}$. We then conclude the convergence of the sequence gives a limit weak solution for the active scalar equation (\ref{S-ase}) in the sense of Definition \ref{def}. It also reveals the failure of constructing non-trivial $L^2$ weak solutions for the stationary even active scalar equation using current scheme. In Section \ref{sec-force}, we present the sum-difference framework of convex integration for the forced equation (\ref{S-ase-g}) and provide heuristics to show that the regularity of convex integration weak solutions does not improve the claimed regularity in Theorem \ref{thm} for the unforced equation. We omit the detail of a complete proof in the forced case, as the proof would be a slight modification of the proof for the forced stationary SQG equation in \cite{DP}.
 


\bigskip

\section{Construction for the unforced active scalar equation}
\label{sec-no-force}

In this section we construct approximating solutions for the active scalar equation (\ref{S-ase}). 
Denote $f=\Gamma^{-1}[\theta]$ and hence $\theta=\Gamma f$. Equation (\ref{S-ase}) with $g\equiv 0$ can be written as
\begin{equation}\label{S-ase2}
\nabla\cdot(\nabla^\perp f\Gamma f)+\Lambda^\gamma \Gamma f=0
\end{equation}
and further in the form
\begin{equation}\label{S-ase3}
\nabla^\perp f\Gamma f-\Lambda^{\gamma-2}\nabla \Gamma f=\nabla^{\perp} V
\end{equation}
for some vector field $V$. The relaxed equation of \eqref{S-ase3} is 
\begin{equation}\label{S-q}
\nabla^\perp f_q\Gamma f_q-\Lambda^{\gamma-2}\nabla \Gamma f_q=\nabla G_q+\nabla^{\perp} V_q, \ \ q\in \mathbb N.
\end{equation}
The iteration process relies on constructing an appropriate increment $W_{q+1}=f_{q+1}-f_q$ to obtain a new solution $f_{q+1}$ associated with a new stress field $G_{q+1}$. That is, $(f_{q+1}, G_{q+1})$ satisfies
\begin{equation}\label{S-q1}
\nabla^\perp f_{q+1}\Gamma f_{q+1}-\Lambda^{\gamma-2}\nabla \Gamma f_{q+1}=\nabla G_{q+1}+\nabla^{\perp} V_{q+1}
\end{equation}
for a new vector $V_{q+1}$ which does not play a role after taking $\div$ on \eqref{S-q1}.
It is not difficult to verify
\begin{equation}\label{G-q1}
\begin{split}
\nabla G_{q+1}=&\left(\nabla^\perp f_{q}\Gamma W_{q+1}+\nabla^\perp W_{q+1}\Gamma f_{q}\right)-\Lambda^{\gamma-2}\nabla \Gamma W_{q+1}\\
&+\left(\nabla^\perp W_{q+1}\Gamma W_{q+1}+\nabla G_q\right)\\
=&: \nabla G_N+ \nabla G_D+\nabla G_O
\end{split}
\end{equation}
where $G_N$, $G_D$ and $G_O$ stand for the Nash error, dissipation error and oscillation error, respectively. 
The essential goal of constructing $W_{q+1}$ is to ensure $G_{q+1}$ is smaller than $G_q$ for all $q\geq 0$, which can be achieved if 
\[\nabla^\perp W_{q+1}\Gamma W_{q+1}+\nabla G_q\]
is small. In other words, the purpose of $W_{q+1}$ is to reduce the error $G_q$.

\subsection{Main iteration}

Fixe $\lambda_0\gg1$. For $b>1$ and $0<\beta<1$, denote 
\[\lambda_{q}=\lceil \lambda_0^{b^q}\rceil, \ \ \delta_q=\lambda_q^{-\beta}, \ \ r_q=(\lambda_q\lambda_{q+1})^{\frac12}.\]
Inductively we expect that $f_q$ and $G_q$ are localized to frequency $\sim \lambda_q$, and the size of $G_q$ is $\delta_q$, i.e. $|G_q|\sim \delta_q$.


We first consider the single plane wave $f(x)=a(x)\cos(\lambda\xi\cdot x)$. We can write
\begin{equation}\notag
\begin{split}
\Gamma f=&\ \frac12\Gamma[ae^{i\lambda\xi\cdot x}]+\frac12\Gamma[ae^{-i\lambda\xi\cdot x}]\\
=&\ \frac12\Gamma_{+}[a]e^{i\lambda\xi\cdot x}+\frac12\Gamma_{-}[a]e^{-i\lambda\xi\cdot x}
\end{split}
\end{equation}
with 
\[\widehat{\Gamma_+[a]}=m(k+\lambda\xi)\widehat a, \ \ \widehat{\Gamma_-[a]}=m(k-\lambda\xi)\widehat a.\]
We further rearrange the terms in $\Gamma f$ to have
\begin{equation}\notag
\begin{split}
\Gamma f=&\ \frac12m(\lambda\xi)a(x)e^{i\lambda\xi\cdot x}+ \frac12m(-\lambda\xi)a(x)e^{-i\lambda\xi\cdot x}\\
&+\frac12\left(\Gamma_{+}[a]- m(\lambda\xi)a(x)\right) e^{i\lambda\xi\cdot x}+\frac12\left(\Gamma_{-}[a]- m(-\lambda\xi)a(x)\right) e^{-i\lambda\xi\cdot x}\\
=&\  im(\lambda\xi)a(x)\sin(\lambda\xi\cdot x)+T_{1,\lambda\xi}[a]\cos(\lambda\xi\cdot x)+T_{2,\lambda\xi}[a]\sin(\lambda\xi\cdot x)
\end{split}
\end{equation}
where we used the odd property of $m$, 
with 
\begin{equation}\notag
\begin{split}
\widehat{T_{1,\lambda\xi}[a]}=&\left( m(k+\lambda\xi)+m(k-\lambda\xi)\right)\widehat a,\\
\widehat{T_{2,\lambda\xi}[a]}=&\left( m(k+\lambda\xi)-m(k-\lambda\xi)-2m(\lambda\xi)\right)\widehat a.
\end{split}
\end{equation}
We observe $\widehat{T_{1,\lambda\xi}[a]}\sim 2m(\lambda\xi)\widehat a$ up to leading order term, which is the symbol for $2\Gamma a$ approximately; while $\widehat{T_{2,\lambda\xi}[a]}$ represents a small error term, since $m$ is odd.

On the other hand we have
\begin{equation}\notag
\nabla^{\perp} f=-a(x)\lambda\xi^{\perp}\sin(\lambda\xi\cdot x)+\nabla^{\perp} a(x)\cos(\lambda\xi\cdot x).
\end{equation}
Therefore,
\begin{equation}\notag
\begin{split}
\nabla^{\perp} f\Gamma f=&-\frac12i\lambda a^2(x)m(\lambda\xi)\xi^{\perp}\\
&+\frac12i\lambda a^2(x)m(\lambda\xi)\xi^{\perp}\cos(2\lambda\xi\cdot x)+\frac12 im(\lambda\xi)a(x)\nabla^{\perp} a(x)\sin(2\lambda\xi\cdot x)\\
&+\left(-\frac12\lambda\xi^{\perp}a(x)T_{1,\lambda\xi}[a]+\frac12\nabla^{\perp} a(x) T_{2,\lambda\xi}[a]\right)\sin(2\lambda\xi\cdot x)\\
&+\nabla^{\perp} a(x)T_{1,\lambda\xi}[a]\cos^2(\lambda\xi\cdot x)-\lambda\xi^{\perp}a(x)T_{2,\lambda\xi}[a]\sin^2(\lambda\xi\cdot x).
\end{split}
\end{equation}
The first term on the right hand side is a low frequency term which is the leading order term and is the source to cancel the aforementioned error $G_q$.

To realize the cancellation, we need the following decomposition lemma.
\begin{Lemma}\label{le-decomp}
Let $G\in C_0^\infty (\mathbb T^2)$ be a vector field. Take $\xi_1=(\frac35, \frac45)$ and $\xi_2=(1,0)$. There exist operators $\mathcal R_1$ and $\mathcal R_2$ of degree 1 such that the decomposition 
\begin{equation}\label{eq-decomp}
\sum_{j=1}^2im(\xi_j)\xi_j^{\perp} \mathcal R_jG=\nabla G+\nabla^{\perp}V
\end{equation}
holds for some vector field $V$.
\end{Lemma}
\pf
Let $\mathcal R_1$ and $\mathcal R_2$ be the operators with Fourier symbols 
\begin{equation}\notag
\widehat { \mathcal R_1}(k)=Ak_1+Bk_2, \ \ \widehat { \mathcal R_2}(k)=Ck_1+Dk_2
\end{equation}
with coefficients $A,B,C, D$ to be determined. 
Acting $\nabla\cdot$ on the identity \eqref{eq-decomp} and taking Fourier transform, we need to solve
\begin{equation}\notag
\begin{split}
\sum_{j=1}^2i^2m(\xi_j)\xi_j^{\perp}\cdot k\widehat{ \mathcal R_j}\widehat G=-|k|^2\widehat G
\end{split}
\end{equation}
which is equivalent to
\begin{equation}\notag
m(\xi_1)(-\frac45k_1+\frac35k_2) \widehat{ \mathcal R_1}+m(\xi_2)k_2 \widehat{ \mathcal R_2}=|k|^2.
\end{equation}
It is easy to find that the following conditions 
\begin{equation}
\begin{cases}
-4Am(\xi_1)=5,\\
3Bm(\xi_1)+5Dm(\xi_2)=5,\\
(3A-4B)m(\xi_1)+5Cm(\xi_2)=0
\end{cases}
\end{equation}
guarantee the identity is satisfied. We note for $m(\xi_1)\neq 0$ and $m(\xi_2)\neq 0$, there are more than one solutions for the coefficients. For instance we can choose
\[A=-\frac{5}{4m(\xi_1)}, \ \ B=\frac{5}{3m(\xi_1)}, \ \ C=\frac{25}{12m(\xi_2)}, \ \ D=\frac{4}{5m(\xi_2)}.\]
It concludes the proof.

\cbdu

Let $X$ be the space equipped with the norm 
\[\|G\|_{X}=\|G\|_{L^\infty}+\sum_{j=1}^2\|\mathcal R_j G\|_{L^\infty}.\]
Denote $\mathbb P_{\approx \lambda}$ by the standard Littlewood-Paley projection operator localized around the frequency $\lambda$. 
The main iteration statement is given below.
\begin{Proposition}\label{prop-1}
Let 
$\lambda_0\gg1$, $0<\gamma<2-\alpha$ and $ \alpha< 1$.
There exist $b>1$ and $0<\beta<1$ satisfying
\begin{equation}\label{parameter}
2b(\alpha-1)+1<\beta< \min\left\{ \frac{2b(2-\gamma)-1}{2b-1}, \frac{2b+2\alpha-3}{2b-1} \right\}
\end{equation}
such that, if $(f_q, G_q, V_q)$ satisfies (\ref{S-q}) with $f_q\in C^\alpha$ and
\begin{equation}\label{iter-1}
f_q=\mathbb P_{\leq 6\lambda_n}f_q,  \ \ G_q=\mathbb P_{\leq 12\lambda_n}G_q, 
\end{equation}
\begin{equation}\label{iter-2}
\|G_q\|_{X}\leq \delta_q,
\end{equation}
\begin{equation}\label{iter-3}
\|G_q\|_{C^{s}}\lesssim \lambda_q^s \delta_q, \ \ \ s\geq \beta,
\end{equation}
there exits $(f_{q+1}, G_{q+1}, V_{q+1})$ satisfying (\ref{S-q}) with $f_{q+1}\in C^\alpha$, and (\ref{iter-1})-(\ref{iter-3}) satisfied with $q$ replaced by $q+1$.
\end{Proposition}
\pf
We first consider the increment
\begin{equation}\notag
W_{q+1}(x)=\sum_{j=1}^2a_{j,q+1}(x)\cos(5\lambda_{q+1}\xi_j\cdot x)
\end{equation}
to motivate the choice of the coefficient functions $a_{j,q+1}(x)$ in the following. Similar analysis as before,
we have
\begin{equation}\notag
\begin{split}
\Gamma W_{q+1}=&\sum_{j=1}^2im(5\lambda_{q+1}\xi_j) a_{j,q+1}\sin(5\lambda_{q+1}\xi_j\cdot x)
+\sum_{j=1}^2T_{1,5\lambda_{q+1}\xi_j}[a_{j,q+1}]\cos(5\lambda_{q+1}\xi_j\cdot x)\\
&+\sum_{j=1}^2T_{2,5\lambda_{q+1}\xi_j}[a_{j,q+1}]\sin(5\lambda_{q+1}\xi_j\cdot x)
\end{split}
\end{equation}
and 
\begin{equation}\notag
\begin{split}
\nabla^{\perp} W_{q+1}=&-\sum_{j=1}^25\lambda_{q+1}\xi_j^{\perp} a_{j,q+1}\sin(5\lambda_{q+1}\xi_j\cdot x)
+\sum_{j=1}^2\nabla^{\perp}a_{j,q+1}\cos(5\lambda_{q+1}\xi_j\cdot x).
\end{split}
\end{equation}
It follows that
\begin{equation}\notag
\begin{split}
\nabla^{\perp} W_{q+1}\Gamma W_{q+1}=&-\sum_{j=1}^225i\lambda_{q+1}^2m(\xi_j)\xi_j^{\perp}a_{j,q+1}^2 \sin^2(5\lambda_{q+1}\xi_j\cdot x)\\
&-\sum_{j\neq j'}25i\lambda_{q+1}^2m(\xi_j)\xi_{j'}^{\perp}a_{j,q+1}a_{j',q+1} \sin(5\lambda_{q+1}\xi_j\cdot x) \sin(5\lambda_{q+1}\xi_{j'}\cdot x)\\
&+\sum_{j,j'}5i\lambda_{q+1}m(\xi_j) a_{j,q+1}\nabla^{\perp}a_{j',q+1}\sin(5\lambda_{q+1}\xi\cdot x)\cos(5\lambda_{q+1}\xi_{j'}\cdot x)\\
&-\sum_{j, j'}5\lambda_{q+1}T_{1,5\lambda_{q+1}\xi_j}[a_{j,q+1}]\xi_{j'}^{\perp}a_{j', q+1}\cos(5\lambda_{q+1}\xi_j\cdot x)\sin(5\lambda_{q+1}\xi_{j'}\cdot x)\\
&+\sum_{j,j'}T_{1,5\lambda_{q+1}\xi_j}[a_{j,q+1}]\nabla^{\perp}a_{j', q+1}\cos(5\lambda_{q+1}\xi_j\cdot x)\cos(5\lambda_{q+1}\xi_{j'}\cdot x)\\
&-\sum_{j,j'}5\lambda_{q+1}T_{2,5\lambda_{q+1}\xi_j}[a_{j,q+1}]\xi_{j'}^{\perp}a_{j',q+1}\sin(5\lambda_{q+1}\xi_j\cdot x)\sin(5\lambda_{q+1}\xi_{j'}\cdot x)\\
&+\sum_{j,j'}T_{2,5\lambda_{q+1}\xi_j}[a_{j,q+1}]\sin(5\lambda_{q+1}\xi_j\cdot x)\nabla^{\perp}a_{j', q+1}\cos(5\lambda_{q+1}\xi_{j'}\cdot x)\\
=:&\ \widetilde J_1+\nabla J_2+\nabla J_3+\nabla J_4+\nabla J_5+\nabla J_6+\nabla J_7.
\end{split}
\end{equation}
The first term $\widetilde J_1$ can be further written as
\begin{equation}\notag
\begin{split}
\widetilde J_1=&-\sum_{j=1}^2\frac{25}{2}i\lambda_{q+1}^2m(\xi_j)\xi_j^{\perp}a_{j,q+1}^2+\sum_{j=1}^2\frac{25}{2}i\lambda_{q+1}^2m(\xi_j)\xi_j^{\perp}a_{j,q+1}^2 \cos(10\lambda_{q+1}\xi_j\cdot x)\\
=&: \widetilde J_{1,1}+\nabla J_{1,2}
\end{split}
\end{equation}
where the low frequency term $\widetilde J_{1,1}$ will be used to cancel the error $G_q$ from the previous level. Indeed, in view of Lemma \ref{le-decomp}, we choose
\begin{equation}\label{a-choice}
a_{j,q+1}=\frac{\sqrt 2}{5}\lambda_{q+1}^{-1}\lambda_q^{\frac12}\delta_q^{\frac12}\left(\mathcal R_j \frac{G_q}{\lambda_q\delta_q}+c_0\right)^{\frac12}
\end{equation}
for a constant $c_0$ such that $\mathcal R_j G_q+\lambda_q\delta_qc_0>0$ for $j=1,2$. We observe 
\[|a_{j,q+1}|\sim \lambda_{q+1}^{-1}\lambda_q^{\frac12}\delta_q^{\frac12}\]
since $G_q$ is localized to frequency $\sim \lambda_q$ with size $\delta_q$ and $\mathcal R_j$ is of order 1.
To ensure $W_{q+1}\in C^{\alpha}$, we impose the condition
\[\lambda_{q+1}^{\alpha} \lambda_{q+1}^{-1}\lambda_q^{\frac12}\delta_q^{\frac12}\lesssim 1.\]
It thus follows 
\begin{equation}\label{alpha}
\alpha<1-\frac1{2b}+\frac{\beta}{2b}.
\end{equation}

To make sure $W_{q+1}$ is localized to frequency $\approx \lambda_{q+1}$, the final form of the increment is
\begin{equation}\notag
W_{q+1}(x)=\sum_{j=1}^2\mathbb P_{\approx \lambda_{q+1}}\left[a_{j,q+1}(x)\cos(5\lambda_{q+1}\xi_j\cdot x)\right]. 
\end{equation}
In view of \eqref{a-choice}, we have
\begin{equation}\label{W-q1}
\|W_{q+1}\|_{X}\lesssim \lambda_{q+1}^{-1}\lambda_q^{\frac12}\delta_q^{\frac12}.
\end{equation}

We proceed to the estimates of the errors $G_N, G_D$ and $G_O$ in \eqref{G-q1}. Note 
\[\nabla G_O=(\widetilde J_{1,1}+\nabla G_q)+\nabla J_{1,2}+\nabla J_2+\nabla J_3+\nabla J_4+\nabla J_5+\nabla J_6+\nabla J_7.\]
We also observe that the items $\nabla G_N, \nabla G_D$ and $\nabla G_O$ have zero mean. 

To estimate $G_N$, applying the inverse operator of the gradient, we obtain
\begin{equation}\notag
G_N=\Delta^{-1}\nabla \cdot \left(\nabla^{\perp}f_q\Gamma W_{q+1}+\nabla^{\perp}W_{q+1}\Gamma f_q\right)
\end{equation}
where the right hand side is localized to $\approx \lambda_{q+1}$. Hence, since $\Gamma$ is of order 1, it follows from \eqref{W-q1}
\begin{equation}\notag
\begin{split}
\|G_N\|_{X}\lesssim &\ \|W_{q+1}\|_{L^\infty}\left(\|\nabla^{\perp}f_{q}\|_{L^\infty}+\|\Gamma f_q\|_{L^\infty}\right)\\
\lesssim& \  \lambda_{q+1}^{-1}\lambda_q^{\frac12}\delta_q^{\frac12}\lambda_q^{1-\alpha}
\end{split}
\end{equation}
where we used $f_{q}\in C^\alpha$ and $\|\nabla^{\perp}f_{q}\|_{C^0}+\|\Gamma f_q\|_{C^0}\lesssim \lambda_q^{1-\alpha}$.

Regarding $G_D$ we have
\begin{equation}\notag
G_D=-\Lambda^{\gamma-2}\Gamma W_{q+1},
\end{equation}
and 
\begin{equation}\notag
\|G_D\|_{X}\lesssim \lambda_{q+1}^{\gamma-1} \|W_{q+1}\|_{X} \lesssim \lambda_{q+1}^{\gamma-1} \lambda_{q+1}^{-1}\lambda_q^{\frac12}\delta_q^{\frac12}\lesssim \lambda_{q+1}^{\gamma-2} \lambda_q^{\frac12}\delta_q^{\frac12}.
\end{equation}

The estimate of $G_O$ takes more effort. We start by applying the decomposition Lemma \ref{le-decomp} to conclude
\begin{equation}\notag
\widetilde J_{1,1} +\nabla G_q=\lambda_{q+1}^2\sum_{j} m(\xi_j)\xi_j^{\perp} \mathbb P_{\leq r_{q+1}} \left(-2a_{j, q+1}\mathbb P_{>r_{q+1}}a_{j,q+1}+ (\mathbb P_{>r_{q+1}}a_{j,q+1})^2\right).
\end{equation}
Hence, by the definition of $a_{j,q+1}$ we deduce 
\begin{equation}\notag
\begin{split}
&\| \Delta^{-1}\nabla\cdot (\widetilde J_{1,1} +\nabla G_q) \|_{X}\\
\lesssim&\ \lambda_{q+1}^{-1+2} (\log r_{q+1})\left( \|a_{j,q+1}\|_{L^\infty}\|\mathbb P_{>r_{q+1}}a_{j,q+1}\|_{L^\infty}+\|\mathbb P_{>r_{q+1}}a_{j,q+1}\|_{L^\infty}^2\right)\\
\lesssim&\ \lambda_{q+1}^{-1+2} (\log r_{q+1})\left( \|a_{j,q+1}\|_{L^\infty}r_{q+1}^{-2}\|\Delta a_{j,q+1}\|_{L^\infty}+r_{q+1}^{-4}\|\Delta a_{j,q+1}\|_{L^\infty}^2\right)\\
\lesssim& \ \lambda_{q+1} (\log r_{q+1})\lambda_{q+1}^{-1}\lambda_q^{\frac12}\delta_q^{\frac12} r_{q+1}^{-2}\lambda_q^2 \lambda_{q+1}^{-1}\lambda_q^{\frac12}\delta_q^{\frac12}\\
\lesssim &\ \lambda_{q+1}^{-1}\lambda_q^3\delta_q r_{q+1}^{-2} \log r_{q+1}.
\end{split}
\end{equation}
Similarly as before, we note
\begin{equation}\notag
J_{1,2}=\Delta^{-1}\nabla\cdot \sum_{j}j\lambda_{q+1}^2m(\xi_j)\xi_j^{\perp} a_{j,q+1}^2,
\end{equation}
and 
\begin{equation}\notag
\|J_{1,2}\|_{X}\lesssim \lambda_{q+1}^{-1+2}\|a_{q+1}\|_{C^0}^2\lesssim \lambda_{q+1} \lambda_{q+1}^{-2}\lambda_q\delta_q.
\end{equation}
Moreover we can write
\begin{equation}\notag
J_{2}=\Delta^{-1}\nabla\cdot \sum_{j\neq j'}25i\lambda_{q+1}^2m(\xi_j)\xi_{j'}^{\perp}a_{j,q+1}a_{j', q+1} \sin(5\lambda_{q+1}\xi_j\cdot x) \sin(5\lambda_{q+1}\xi_{j'}\cdot x),
\end{equation}
\begin{equation}\notag
J_{3}=\Delta^{-1}\nabla\cdot \sum_{j\neq j'}5i\lambda_{q+1}m(\xi_j)a_{j,q+1}\nabla^{\perp}a_{j', q+1} \sin(5\lambda_{q+1}\xi_j\cdot x) \cos(5\lambda_{q+1}\xi_{j'}\cdot x),
\end{equation}
and thus
\begin{equation}\notag
\|J_{2}\|_{X}\lesssim \lambda_{q+1}^{-1+2}\|a_{q+1}\|_{C^0}^2\lesssim \lambda_{q+1} \lambda_{q+1}^{-2}\lambda_q\delta_q,
\end{equation}
\begin{equation}\notag
\|J_{3}\|_{X}\lesssim \lambda_{q+1}^{-1+1}\lambda_q\|a_{q+1}\|_{C^0}^2\lesssim \lambda_{q} \lambda_{q+1}^{-2}\lambda_q\delta_q.
\end{equation}
As observed earlier, the Fourier symbol of $T_{1,\lambda\xi}[a]$ satisfies $\widehat {T_{1,\lambda\xi}[a]}\sim 2m(\lambda\xi)\widehat a$. We thus have
\begin{equation}\notag
\|T_{1, 5\lambda_{q+1}\xi_j}[a_{j,q+1}]\|_X\lesssim |m(5\lambda_{q+1}\xi_j)|\|a_{j,q+1}\|_{X}\lesssim \lambda_{q+1} |m(\xi_j)| \|a_{j,q+1}\|_{L^\infty},
\end{equation}
which implies
\begin{equation}\notag
\begin{split}
\|J_4\|_{X}\lesssim&\ \lambda_{q+1}^{-1+1}\lambda_{q+1} \|a_{j,q+1}\|_{L^\infty}^2\lesssim \lambda_{q+1} \lambda_{q+1}^{-2}\lambda_q\delta_q\lesssim  \lambda_{q+1}^{-1}\lambda_q\delta_q,\\
\|J_5\|_{X}\lesssim&\ \lambda_{q+1}^{-1}\lambda_q\lambda_{q+1} \|a_{j,q+1}\|_{L^\infty}^2\lesssim \lambda_{q} \lambda_{q+1}^{-2}\lambda_q\delta_q\lesssim  \lambda_{q+1}^{-2}\lambda_q^2\delta_q.
\end{split}
\end{equation}
We also observe that $T_{2, \cdot}$ is a minor error term compared to $T_{1, \cdot}$. Hence we have
\[\|J_6\|_{X}\lesssim \|J_4\|_{X}\lesssim  \lambda_{q+1}^{-1}\lambda_q\delta_q, \]
\[\|J_7\|_{X}\lesssim \|J_5\|_{X}\lesssim  \lambda_{q+1}^{-2}\lambda_q^2\delta_q. \]

Summarizing the estimates above gives the estimate for the new error
\begin{equation}\notag
\|G_{q+1}\|_X\lesssim \lambda_{q+1}^{-1}\lambda_q^{\frac12}\delta_q^{\frac12}\lambda_q^{1-\alpha}+\lambda_{q+1}^{\gamma-2}\lambda_q^{\frac12}\delta_q^{\frac12}+\lambda_{q+1}^{-1}\lambda_q^3\delta_q r_{q+1}^{-2}\log r_{q+1}+\lambda_{q+1}^{-1}\lambda_q\delta_q.
\end{equation}
Therefore, to show the estimate \eqref{iter-2} holds for $q+1$, we only need to verify 
\begin{equation}\notag
\begin{cases}
\lambda_{q+1}^{-1}\lambda_q^{\frac12}\delta_q^{\frac12}\lambda_q^{1-\alpha}\lesssim \delta_{q+1}\\
\lambda_{q+1}^{\gamma-2}\lambda_q^{\frac12}\delta_q^{\frac12} \lesssim \delta_{q+1}\\
\lambda_{q+1}^{-1}\lambda_q^3\delta_q r_{q+1}^{-2}\log r_{q+1}\lesssim \delta_{q+1}\\
\lambda_{q+1}^{-1}\lambda_q\delta_q \lesssim \delta_{q+1}
\end{cases}
\end{equation}
which is equivalent to, (recalling $r_{q+1}=(\lambda_q\lambda_{q+1})^{\frac12}$)
\begin{equation}\label{para-ss}
\begin{cases}
b\beta-b+\frac12-\frac12\beta+1-\alpha<0\\
b\beta+b(\gamma-2)+\frac12-\frac12\beta <0\\
b\beta-2b+2-\beta<0\\
b\beta-b+1-\beta<0.
\end{cases}
\end{equation}
Since $\alpha<1-\frac{1}{2b}+\frac{\beta}{2b}$ from \eqref{alpha}, it follows from the first inequality of \eqref{para-ss} that 
\begin{equation}\notag
\beta<1+\frac{\alpha-1}{b-\frac12}< 1+ \frac{-\frac{1}{2b}+\frac{\beta}{2b}}{b-\frac12}
\end{equation}
which is equivalent to
\begin{equation}\notag
(2b+1)(b-1)(\beta-1)<0.
\end{equation}
It is satisfied for all $b>1$ and $0<\beta<1$. Hence when $b=1^+$, it yields $\alpha<1-\frac12+\frac12<1$. The second inequality of (\ref{para-ss}) gives
\begin{equation}\notag
\gamma<2-\beta-\frac1{2b}+\frac{\beta}{2b}<\frac32-\frac12\beta
\end{equation}
for $b=1^+$. The third and forth inequalities of (\ref{para-ss}) yield 
\begin{equation}\notag
(\beta-2)(b-1)<0, \ \ (\beta-1)(b-1)<0
\end{equation}
which are satisfied for all $b>1$ and $0<\beta<1$. On the other hand, we notice the conditions on the parameters of the proposition are compatible with \eqref{para-ss}. Hence the estimate \eqref{iter-2} for $q+1$ is verified. 

The claim $f_{q+1}\in C^{\alpha}$ follows from \eqref{alpha} and \eqref{W-q1}. The frequency support property of \eqref{iter-1} for $q+1$ follows directly from the definition of $W_{q+1}$. While the estimate \eqref{iter-3} for $q+1$ follows from \eqref{iter-2} with $q+1$ and \eqref{iter-1} with $q+1$.

\cbdu

\medskip

\subsection{Proof of Theorem \ref{thm}} 
Under the assumptions on the parameters as in Proposition \ref{prop-1}, we first choose $f_0=0$ and $G_0=0$ which satisfies equation \eqref{S-q} with a zero-vector field $V_0$. Note $f_0$ and $G_0$ also satisfy \eqref{iter-1}-\eqref{iter-3}. Applying Proposition \ref{prop-1} inductively we obtain a sequence of solutions $\{f_q, G_q, V_q\}$ to \eqref{S-q} satisfying \eqref{iter-1}-\eqref{iter-3} and $f_q\in C^\alpha$. Taking the limit $q\to \infty$, it is clear that $G_q\to 0$ in $L^1$. 
We need to show that as $q\to \infty$, there is a limit function $f$ of the sequence $\{f_q\}$ such that $\theta=\Gamma f$ is a weak solution of \eqref{S-ase} in the sense of Definition \ref{def}. Indeed, for any $\psi\in C^\infty(\mathbb T^2)$ with $\psi=\mathbb P_{\leq \Lambda_0}\psi$ for some frequency number $\Lambda_0$, we have
\begin{equation}\notag
\int_{\mathbb T^2}\theta_q u_q\cdot \nabla \psi \, dx=\int_{\mathbb T^2}\mathbb P_{\leq \Lambda_0}(\theta_q u_q)\cdot \nabla \psi \, dx\to \int_{\mathbb T^2}\mathbb P_{\leq \Lambda_0}(\theta T[\theta])\cdot \nabla \psi \, dx \ \ \mbox{as} \ \ q\to \infty,
\end{equation}
 due to the property of compact support in frequency space.

\bigskip

\section{The forced stationary equation}
\label{sec-force}

In this section we move on to construct non-unique weak solutions for the forced stationary active scalar equation 
\begin{equation}\label{S-ase-g}
\begin{split}
u \cdot\nabla \theta+\Lambda^\gamma \theta=&\ g, \\
u=T[\theta]=&\ \nabla^{\perp}\Gamma^{-1}[\theta].
\end{split}
\end{equation}
We first describe the convex integration scheme applied to the sum-difference formation of the forced equation, which was previously exploited in \cite{DF, DP, DP2}. Then we state the main iteration proposition, which will lead to the conclusion of existence of non-unique weak solutions with regularity consistent as stated in the proposition. The purpose is to show that, peculiarly, the presence of a forcing does not necessarily improve the regularity of convex integration weak solutions.

\subsection{Sum-difference formulation}
Assume $(\theta, g_1)$ and $(\widetilde\theta, g_2)$ are two solutions of \eqref{S-ase-g}. 
Denote
\[\eta=\Gamma^{-1}\theta, \ \ \widetilde\eta=\Gamma^{-1}\widetilde\theta,\ \ \Pi=\frac12(\eta+\widetilde \eta), \ \ \mu=\frac12(\eta-\widetilde \eta),\] 
and
\[p=\Gamma \Pi, \ \ m=\Gamma \mu, \ \ T[p]=\nabla^{\perp}\Pi, \ \ T[m]=\nabla^{\perp}\mu.\]
It follows from \eqref{S-ase-g} that
\begin{equation}\label{pm4}
\begin{split}
\Gamma \Pi\nabla^{\perp}\Pi+\Gamma \mu\nabla^{\perp}\mu- \Lambda^{\gamma-2} \nabla\Gamma\Pi=&\ \nabla  G+\nabla^{\perp}V,\\
\Gamma \mu\nabla^{\perp}\Pi+\Gamma \Pi\nabla^{\perp}\mu-\Lambda^{\gamma-2}\nabla \Gamma\mu=&\ \nabla\widetilde G+\nabla^{\perp}\widetilde V
\end{split}
\end{equation} 
with vector fields $G$ and $\widetilde G$ satisfying $\Delta G=\frac12(g_1+g_2)$ and $\Delta \widetilde G=\frac12(g_1-g_2)$, and some vector fields $V$ and $\widetilde V$. 
The relaxed system of \eqref{pm4} is 
\begin{equation}\label{pm-q}
\begin{split}
\Gamma \Pi_q\nabla^{\perp}\Pi_q+\Gamma \mu_q\nabla^{\perp}\mu_q- \Lambda^{\gamma-2} \nabla\Gamma\Pi_q=&\ \nabla G_q+\nabla^{\perp} V_q,\\
\Gamma \mu_q\nabla^{\perp}\Pi_q+\Gamma \Pi_q\nabla^{\perp}\mu_q-\Lambda^{\gamma-2}\nabla\Gamma \mu_q=&\ \nabla \widetilde G_q+\nabla^{\perp}\widetilde V_q.
\end{split}
\end{equation}

To obtain a sequence of approximating solutions $\{(\Pi_q,\mu_q, G_q, \widetilde G_q)\}_{q\geq 0}$ of (\ref{pm-q}), we will perform an active convex integration scheme to the second equation and a passive scheme for the first equation, which will be described in detail below. The goal is to have $\widetilde G_q$ approaches zero as $q\to\infty$. Consequently the limit $(\Pi, \mu, G, 0)$ with non-vanishing $\mu$ solves (\ref{pm4}). It then follows that there are two distinct stationary solutions $\theta=\Gamma(\Pi+\mu)$ and $\widetilde\theta=\Gamma(\Pi-\mu)$ of (\ref{S-ase}) with forcing $g=\Delta G$.

Denote
\[\theta_q=\Gamma (\Pi_q+\mu_q), \ \ \widetilde\theta_q=\Gamma (\Pi_q-\mu_q), \ \ g_{1,q}=\Delta(G_q+\widetilde G_q), \ \ g_{2,q}=\Delta(G_q-\widetilde G_q). \]
It is obvious that $(\theta_q, g_{1,q})$ and $(\widetilde\theta_q, g_{2,q})$ satisfy the stationary forced active scalar equation (\ref{S-ase-g}). 
We describe the construction from $q$-th to $(q+1)$-th level and from $(q+1)$-th to $(q+2)$-th level in the following.
 
We will design perturbations $M_{q+1}$ and $M_{q+2}$ such that
\[\mu_{q+1}=\mu_q+M_{q+1}, \ \ \Pi_{q+1}=\Pi_q-M_{q+1},\]
\[\mu_{q+2}=\mu_{q+1}+M_{q+2}, \ \ \Pi_{q+2}=\Pi_{q+1}+M_{q+2}.\]
Straightforward computation shows that
\begin{equation}\notag
\begin{split}
\Gamma^{-1}\theta_{q+1}=&\ \eta_{q+1}= \Pi_{q+1}+\mu_{q+1}=\Pi_q+\mu_q\\
=&\ \eta_q=\Lambda^{-1}\theta_{q}, \\
\Gamma^{-1}\widetilde\theta_{q+1}=&\ \widetilde\eta_{q+1}= \Pi_{q+1}-\mu_{q+1}=\Pi_q-\mu_q-2M_{q+1}\\
=&\ \widetilde\eta_q-2M_{q+1}=\Lambda^{-1}\widetilde\theta_{q}-2M_{q+1},\\
\Gamma^{-1}\theta_{q+2}=&\ \eta_{q+2}= \Pi_{q+2}+\mu_{q+2}=\Pi_{q+1}+\mu_{q+1}+2M_{q+2}\\
=&\ \eta_{q+1}+2M_{q+2}=\Lambda^{-1}\theta_{q+1}+2M_{q+2}, \\
\Gamma^{-1}\widetilde\theta_{q+2}=&\ \widetilde\eta_{q+2}= \Pi_{q+2}-\mu_{q+2}=\Pi_{q+1}-\mu_{q+1}\\
=&\ \widetilde\eta_{q+1}=\Lambda^{-1}\widetilde\theta_{q+1}.
\end{split}
\end{equation}
Let $G_{q+1}$, $\widetilde G_{q+1}$, $G_{q+2}$ and $\widetilde G_{q+2}$ be the new stress errors respectively for the solutions $\Pi_{q+1}$ and $\mu_{q+1}$, $\Pi_{q+2}$ and $\mu_{q+2}$. 
Applying (\ref{pm-q}) gives 
\begin{equation}\label{g1-q1}
\begin{split}
\nabla G_{q+1}=&-\nu\Lambda^{\gamma-1}\nabla M_{q+1}-\left(\Lambda \widetilde\eta_q\nabla^{\perp} M_{q+1}+ \Lambda M_{q+1}\nabla^{\perp} \widetilde\eta_q\right)\\
&+2\Lambda M_{q+1}\nabla^{\perp}M_{q+1}+\nabla G_q,
\end{split}
\end{equation}
\begin{equation}\label{g1-q2}
\begin{split}
\nabla G_{q+2}=&\ \nu\Lambda^{\gamma-1}\nabla M_{q+2}+\left(\Lambda \eta_{q+1}\nabla^{\perp} M_{q+2}+ \Lambda M_{q+2}\nabla^{\perp} \eta_{q+1}\right)\\
&+2\Lambda M_{q+2}\nabla^{\perp}M_{q+2}+\nabla G_{q+1},
\end{split}
\end{equation}
\begin{equation}\label{g-q1}
\begin{split}
\nabla\widetilde G_{q+1}=&-\Lambda^{\gamma-1}\nabla \Gamma M_{q+1}+\left(\Gamma \widetilde\eta_q\nabla^{\perp} M_{q+1}+ \Gamma M_{q+1}\nabla^{\perp} \widetilde\eta_q\right)\\
&+\left(\nabla\widetilde G_q-2\Gamma M_{q+1}\nabla^{\perp}M_{q+1}\right),
\end{split}
\end{equation}
\begin{equation}\label{g-q2}
\begin{split}
\nabla\widetilde G_{q+2}=&-\Lambda^{\gamma-1}\nabla \Gamma M_{q+2}+\left(\Gamma \eta_{q+1}\nabla^{\perp} M_{q+2}+ \Gamma M_{q+2}\nabla^{\perp} \eta_{q+1}\right)\\
&+\left(\nabla\widetilde G_{q+1}+2\Gamma M_{q+2}\nabla^{\perp}M_{q+2}\right).
\end{split}
\end{equation}

Similarly as in Section \ref{sec-no-force}, $M_{q+1}$ and $M_{q+2}$ will be designed to ensure 
\begin{equation}\label{cancel1}
\nabla\widetilde G_q-2\Gamma M_{q+1}\nabla^{\perp}M_{q+1}\sim \nabla^{\perp} V_{q+1}
\end{equation}
\begin{equation}\label{cancel2}
\nabla\widetilde G_{q+1}+2\Gamma M_{q+2}\nabla^{\perp}M_{q+2}\sim \nabla^{\perp} V_{q+2}
\end{equation}
up to small errors, for some vector fields $V_{q+1}$ and $V_{q+2}$.

The advantage of allowing external forcing relies on the property of this particular scheme
\begin{equation}\label{skip}
\eta_{q+1}=\eta_q, \ \ \widetilde \eta_{q+2}=\widetilde \eta_{q+1}, \ \ \mbox{for any even} \ \ q\geq 0.
\end{equation}
This directly improves the estimates for the Nash errors in (\ref{g-q1}) and (\ref{g-q2}) compared to the case in Section \ref{sec-no-force}. 


The main iteration statement is

\begin{Proposition}\label{prop-2}
Assume $\lambda_0\gg1$, $0<\gamma<2-\alpha$ and $\alpha<1$.
There exist $b>1$ and $0<\beta<1$ satisfying
\begin{equation}\label{parameter-f}
2b(\alpha-1)+1<\beta< \min\left\{ \frac{2b(2-\gamma)-1}{2b-1}, \frac{2b+2\alpha-3}{2b-1} \right\}
\end{equation}
such that, if $(\Pi_q, \mu_q, G_q, \widetilde G_q)$ satisfies (\ref{pm-q}) with $\Pi_q, \mu_q\in C^\alpha$ and
\begin{equation}\label{iter-1-f}
\Pi_q=\mathbb P_{\leq 6\lambda_q}\Pi_q, \ \ \mu_q=\mathbb P_{\leq 6\lambda_q}\mu_q, \ \ G_q=\mathbb P_{\leq 12\lambda_q}G_q, \ \ \widetilde G_q=\mathbb P_{\leq 12\lambda_q}\widetilde G_n,
\end{equation}
\begin{equation}\label{iter-2-f}
\|G_q\|_{X}\leq 1-\delta_q^{\frac12},
\end{equation}
\begin{equation}\label{iter-3-f}
\|G_q\|_{C^{s}}\lesssim \lambda_q^s \delta_q, \ \ \ s\geq \beta,
\end{equation}
\begin{equation}\label{iter-4-f}
\|\widetilde G_q\|_{X}\leq \delta_q
\end{equation}
 there exits $(\Pi_{q+1}, \mu_{q+1}, G_{q+1}, \widetilde G_{q+1})$ satisfying (\ref{pm-q}) with $\Pi_{q+1}, \mu_{q+1}\in C^\alpha$, and (\ref{iter-1-f})-(\ref{iter-4-f}) satisfied with $q$ replaced by $q+1$.
\end{Proposition}

\medskip

\subsection{Heuristics for the proof of Proposition \ref{prop-2}}
\label{sec-heuristics}
Based on the proof of Proposition \ref{prop-1}, Proposition \ref{prop-2} can be proved in a close analogous way as the Proposition 3.6 from \cite{DP} (or the iteration propositions from \cite{DF, DP2}). To avoid repetition, we only give the precise construction of $M_{n+1}$ and $M_{n+2}$, and provide heuristic analysis to show the estimates claimed in Proposition \ref{prop-2}.

The perturbations $M_{n+1}$ and $M_{n+2}$ have the form
\begin{equation}\label{m-q1-construct}
\begin{split}
M_{q+1}(x)=& \sum_{j=1}^2P_{\leq r_{q+1}} \left(a_{j, q+1}(x)\right)\cos(5\lambda_{q+1}\xi_j\cdot x),\\
M_{q+2}(x)=& \sum_{j=1}^2P_{\leq r_{q+2}} \left(a_{j, q+2}(x)\right)\cos(5\lambda_{q+2}\xi_j\cdot x)
\end{split}
\end{equation}
with the coefficient functions 
\begin{equation}\label{a}
\begin{split}
a_{j, q+1}=&\ \frac{\sqrt 2}{5}\lambda_{q+1}^{-1}\lambda_q^{\frac12}\delta_q^{\frac12}\left(c_0-\mathcal R_j\left(\frac{\widetilde G_q}{\lambda_q\delta_q}\right) \right)^{\frac12},\\
a_{j, q+2}=&\ \frac{\sqrt 2}{5}\lambda_{q+2}^{-1}\lambda_{q+1}^{\frac12}\delta_{q+1}^{\frac12}\left(c_0+\mathcal R_j\left(\frac{\widetilde G_{q+1}}{\lambda_{q+1}\delta_{q+1}}\right) \right)^{\frac12}
\end{split}
\end{equation}
for a constant $c_0\geq 100$ such that the quantities under the square root are positive. 
Such choice of the perturbations will make the cancellations (\ref{cancel1}) and (\ref{cancel2}) achieved, similarly as shown in Section 2.

It is clear that
\begin{equation}\label{m-q1}
|M_{q+1}|\sim \lambda_{q+1}^{-1}\lambda_q^{\frac12}\delta_q^{\frac12}, \ \  |M_{q+2}|\sim \lambda_{q+2}^{-1}\lambda_{q+1}^{\frac12}\delta_{q+1}^{\frac12}.
\end{equation}
For $\alpha$ satisfying 
\begin{equation}\notag
\alpha<1-\frac1{2b}+\frac{\beta}{2b},
\end{equation}
the perturbation $M_{n+1}$ is in $C^{\alpha}$ since  
\[\lambda_{q+1}^{\alpha} \lambda_{q+1}^{-1}\lambda_q^{\frac12}\delta_q^{\frac12}\lesssim \lambda_q^{\alpha b-b+\frac12-\frac12\beta}\lesssim 1.\]

Thanks to \eqref{skip}, the error equation of (\ref{g-q1}) is essentially 
\begin{equation}\label{est-gn1}
\begin{split}
\nabla\widetilde G_{q+1}=&-\Lambda^{\gamma-2}\nabla \Gamma M_{q+1}+\left(\Gamma \widetilde\eta_{q-1}\nabla^{\perp} M_{q+1}+ \Gamma M_{q+1}\nabla^{\perp} \widetilde\eta_{q-1}\right)\\
&+\left(\nabla\widetilde G_q-2\Gamma M_{q+1}\nabla^{\perp}M_{q+1}\right)\\
=:&\ \nabla \widetilde G_D+\nabla \widetilde G_N+\nabla \widetilde G_O.
\end{split}
\end{equation}
We observe that $\widetilde G_D$ and $\widetilde G_O$ can be estimated exactly as $G_D$ and $G_O$ from Section \ref{sec-no-force}. Hence we have
\[ \|\widetilde G_D\|_X\lesssim \lambda_{q+1}^{\gamma-2}\lambda_q^{\frac12}\delta_q^{\frac12}, \]
\[ \|\widetilde G_O\|_X\lesssim \lambda_{q+1}^{-1}\lambda_q^3\delta_q r_{q+1}^{-2} \log r_{q+1}+\lambda_q\lambda_{q+1}^{-1}\delta_q. \]
The essential difference compared to the non-forced case comes from the Nash error $\widetilde G_N$. Indeed, recalling $\widetilde\eta_{q-1}=\Pi_{q-1}-\mu_{q-1}$ and
applying the inductive assumption $\Pi_{q-1}, \mu_{q-1}\in C^\alpha$, we have
\[\|\Gamma \widetilde\eta_{q-1}\|_{L^\infty}\lesssim \|\Gamma\Pi_{q-1}\|_{L^\infty}+\|\Gamma\mu_{q-1}\|_{L^\infty}\lesssim \lambda_{q-1}^{1-\alpha}\] 
and similarly 
\[\|\nabla^{\perp} \widetilde\eta_{q-1}\|_{L^\infty}\lesssim \|\nabla^{\perp}\Pi_{q-1}\|_{L^\infty}+\|\nabla^{\perp}\mu_{q-1}\|_{L^\infty}\lesssim \lambda_{q-1}^{1-\alpha}.\] 
Thus combining with \eqref{m-q1} we obtain
\begin{equation}\notag
\begin{split}
\|\widetilde G_N\|_{X}\lesssim &\ \|M_{n+1}\|_{L^\infty}\left(\|\nabla^{\perp}\widetilde \eta_{n-1}\|_{L^\infty}+\|\Gamma \widetilde\eta_{n-1}\|_{L^\infty} \right)\\
\lesssim& \ \lambda_{q+1}^{-1}\lambda_q^{\frac12}\delta_q^{\frac12}\lambda_{q-1}^{1-\alpha}.
\end{split}
\end{equation}
Therefore we have 
\[\|\widetilde G_{q+1}\|_{X}\lesssim \lambda_{q+1}^{\gamma-2}\lambda_q^{\frac12}\delta_q^{\frac12}+\lambda_{q+1}^{-1}\lambda_q^3\delta_q r_{q+1}^{-2} \log r_{q+1}+\lambda_q\lambda_{q+1}^{-1}\delta_q+\lambda_{q+1}^{-1}\lambda_q^{\frac12}\delta_q^{\frac12}\lambda_{q-1}^{1-\alpha}.\]

Again, to carry on the iteration process, we need to make sure that $\|\widetilde G_{q+1}\|_{X}\lesssim \delta_{q+1}$.
Following a similar analysis for the parameters as in Section \ref{sec-no-force}, we need to impose 
\begin{equation}\notag
\begin{cases}
\lambda_{q+1}^{-1}\lambda_q^{\frac12}\delta_q^{\frac12}\lambda_{q-1}^{1-\alpha}\lesssim \delta_{q+1}\\
\lambda_{q+1}^{\gamma-2}\lambda_q^{\frac12}\delta_q^{\frac12} \lesssim \delta_{q+1}\\
\lambda_{q+1}^{-1}\lambda_q^3\delta_q r_{q+1}^{-2}\log r_{q+1}\lesssim \delta_{q+1}\\
\lambda_{q+1}^{-1}\lambda_q\delta_q \lesssim \delta_{q+1}.
\end{cases}
\end{equation}
By the definitions of $\lambda_q$, $\delta_q$ and $r_q$, the system above is satisfied provided that
\begin{equation}\label{para-ss-f}
\begin{cases}
b\beta-b+\frac12-\frac12\beta+\frac{1}{b}(1-\alpha)<0\\
b\beta+b(\gamma-2)+\frac12-\frac12\beta <0\\
b\beta-2b+2-\beta<0\\
b\beta-b+1-\beta<0.
\end{cases}
\end{equation}
Comparing \eqref{para-ss-f} with \eqref{para-ss}, the difference is in the first inequality. Note that the last inequality of \eqref{para-ss-f} gives the most stringent condition for $\beta$, that is, $\beta<1$. Hence \eqref{para-ss-f} yields essentially the same conditions on parameters as in Section \ref{sec-no-force} for the unforced case. 

Through the analysis above, we do not see an improvement on the regularity of the weak solutions for the forced equation, constructed via the convex integration scheme. The key observation is that, although the estimate for the Nash error can be improved, the estimates for the oscillation errors can not be improved. This discovery exposes the intrinsic difference in the nonlinear structures between the even active scalar equations and other fluid equations, like the odd active scalar equations and the Euler equation.






\bigskip



\begin{thebibliography}{XX}







\bibitem{BDLSV}
T. Buckmaster, C. De Lellis, L. Sz\'ekelyhidi, and V. Vicol.
\newblock {\em Onsager's conjecture for admissible weak solutions}.
\newblock Comm. Pure Appl. Math., https://doi.org/10.1002/cpa.21781. 2018.


\bibitem{BV19}
T. Buckmaster, and V. Vicol.
\newblock {\em Convex integration and phenomenologies in turbulence}.
\newblock Journal of EMS Surveys in Mathematical Sciences, 6 (1-2): 143--263, 2019.


\bibitem{BHP}
A. Bulut, M.K. Huynh, and S. Palasek.
\newblock {\em Convex integration above the Onsager exponent for the forced Euler equations}.
\newblock arXiv:2301.00804, 2023.

\bibitem{BHP2}
A. Bulut, M.K. Huynh, and S. Palasek.
\newblock {\em Non-uniqueness up to the Onsager threshold for the forced SQG equation}.
\newblock arXiv:2310.12947, 2023.





\bibitem{CKL}
X. Cheng, H. Kwon, and D. Li.
\newblock {\em Non-uniqueness of steady-state weak solutions to the surface quasi-geostrophic equations}.
\newblock Commun. Math. Phys. 388, 1281--1295, 2021.










\bibitem{CCF}
A. Castro, D. C\'ordoba, and D. Faraco.
\newblock {\em Mixing solutions for the Muskat problem}.
\newblock Inventiones mathematicae, 226:251--348, 2021.


\bibitem{CFG}
A. C\'ordoba, D. Faraco and F. Gancedo.
\newblock {\em Lack of uniqueness for weak solutions of the incompressible porous media equation}.
\newblock Arch. Ration. Mech. Anal., 200(3):725--746, 2011.

\bibitem{CFM}
A. Castro, D. Faraco, and F. Mengual.
\newblock {\em Localized mixing zone for Muskat bubbles and turned interfaces}.
\newblock Ann. PDE, 8(1), 7, 50, 2022.

\bibitem{CGO}
D. C\'ordoba, F. Gancedo and R. Orive.
\newblock {\em Analytical behavior of two-dimensional incompressible flow in porous media}.
\newblock J. Math. Phys., 48(6), 065206, 2007.

\bibitem{DF}
M. Dai and S. Friedlander.
\newblock {\em Non-uniqueness of forced active scalar equations with even drift operators}.
\newblock arXiv:2311.06064, 2023.

\bibitem{DP}
M. Dai and Q. Peng.
\newblock {\em Non-unique stationary solutions of forced SQG}.
\newblock arXiv:2302.03283, 2023.

\bibitem{DP2}
M. Dai and Q. Peng.
\newblock {\em Non-unique weak solutions of forced SQG}.
\newblock arXiv:2310.13537, 2023.




\bibitem{DLS1}
C. De Lellis, and L. Sz\'ekelyhidi.
\newblock {\em Dissipative continuous Euler flows}.
\newblock Invent. Math., Vol.193 No. 2: 377--407, 2013.

\bibitem{DLS2}
C. De Lellis, and L. Sz\'ekelyhidi.
\newblock {\em The Euler equations as a differential inclusion}.
\newblock Ann. of Math.,  Vol.170 No.3: 1417--1436, 2009.




\bibitem{FV1}
S. Friedlander and V. Vicol.
\newblock {\em Global well-posedness for an advection-diffusion equation arising in magneto-geostrophic dynamics}.
\newblock Ann. Inst. H. Poincar\'e Anal. Non Lin\'eaire, 28(2): 283-301, 2011.


\bibitem{GR}
V. Giri and R. O. Radu.
\newblock {\em The 2D Onsager conjecture: a Newton-Nash iteration}.
\newblock arXiv:2305.18105, 2023.






\bibitem{Is}
P. Isett.
\newblock {\em A Proof of Onsager's Conjecture}.
\newblock Ann. of Math., Vol.188 No.3: 1--93, 2018.



\bibitem{IV}
P. Isett and V. Vicol.
\newblock {\em H\"older continuous solutions of active scalar equations}.
\newblock Ann. PDE., 1(1): 1--77, 2015.









\bibitem{Luo}
X. Luo.
\newblock {\em Stationary solutions and nonuniqueness of weak solutions for the Navier-Stokes equations in high dimensions}.
\newblock Arch. Ration. Mech. Anal., 233:701--747, 2019.








\bibitem{Moff}
H. K. Moffatt.
\newblock {\em Magnetostrophic turbulence and the geodynamo}.
\newblock In: IUTAM Symposium on Computational Physics and New Perspectives in Turbulence, vol. 4 of IUTAM Bookser. Springer, Dordrecht, pp 339-346, 2008.

\bibitem{Nash}
J. Nash. 
\newblock {\em $C^1$ isometric embeddings}.
\newblock Ann. Math., 60:383–396, 1954.










\bibitem{Shv}
R. Shvydkoy.
\newblock {\em Convex integration for a class of active scalar equations}.
\newblock J. Am. Math. Soc., 24(4): 1159--1174, 2011.

\bibitem{Sze}
L. Sz\'ekelyhidi Jr.
\newblock {\em Relaxation of the incompressible porous media equation}.
\newblock Ann. Sci. \'Ec. Norm. Sup\'er., (4) 45(3): 491--509, 2012.






\end{thebibliography}
\end{document}